\documentclass{amsart}

\usepackage{amsfonts}
\usepackage{amssymb}
\usepackage{mathrsfs} %no en windows
\usepackage{amsmath}
\usepackage[english]{babel} %si en windows
\usepackage[latin1]{inputenc}
\usepackage[all]{xy}
\newcommand{\pn}{\par\noindent}

\newcommand{\edge}{\ar@{-}}

\newcommand{\fd}{\mathrm{fin.dim.}\,}
\newcommand{\rad}{\mathrm{rad}\,}
\newcommand{\soc}{\mathrm{soc}\,}
\newcommand{\mad}{\mathrm{mod}\,}

\newcommand{\add}{\mathrm{add}\,}
\newcommand{\pd}{\mathrm{pd}\,}

\newcommand{\tap}{\mathrm{top}\,}
\newcommand{\SM}{S(M)}
\newcommand{\QM}{Q(M)}

\newcommand{\Ker}{\mathrm{Ker}\,}
\newcommand{\Coker}{\mathrm{Coker}\,}
\newcommand{\LL}{\ell\ell\,}
\newcommand{\F}{\mathcal{F}\,}

{}
\newenvironment{dem}{\noindent {\bf Proof.}}{\hfill $\Box$\\}

\newtheorem{pro}{Proposition}[section]
\newtheorem{teo}[pro]{Theorem}
\newtheorem{defi}[pro]{Definition}
\newtheorem{lem}[pro]{Lemma}
\newtheorem{cor}[pro]{Corollary}

\newtheorem{ex}[pro]{Example}
\usepackage[all]{xy}
\newtheorem{reNN}{Remark} {}

\begin{document}
\title{Finitistic dimension through infinite projective dimension.}
\author[F.~Huard - M.~Lanzilotta - O.~Mendoza ]
{Fran\c{c}ois Huard, Marcelo Lanzilotta, Octavio Mendoza}

%%%%%%%%%%%%%%%%%%%%%%%%%%%%%%%%%%%%%%%%%%%%%%%%%%%%%%%%%%%%%%%%

\address{\pn Fran\c cois Huard; Department of mathematics,
Bishop's University. Sherbrooke, Qu\'ebec, Canada,  J1M1Z7.}
\email{fhuard@ubishops.ca}
\address{\pn Marcelo Lanzilotta; Centro de Matem\'atica (CMAT), Igu\'a 4225,
Universidad de la Rep\'ublica. CP 11400, Montevideo, Uruguay.}
\email{marclan@cmat.edu.uy}
\address{\pn Octavio Mendoza Hern\'andez; Instituto de Matem\'aticas,
Universidad Nacionale Aut\'onoma de M\'exico. Circuito Exterior,
Ciudad Universitaria, C.P. 04510, M\'exico, D.F. M\'EXICO.}
\email{omendoza@matem.unam.mx}

%%%%%%%%%%%%%%%%%%%%%%%%%%%%%%%%%%%%%%%%%%%%%%%%%%%%%%%%%%%%%%%%
\begin{abstract}
We show that an Artin algebra $\Lambda$ having at most three radical
layers of infinite projective dimension has finite finitistic
dimension, generalizing the known result for algebras with vanishing
radical cube.

\end{abstract}

\thanks{The last author was partially supported by the project PAPIIT-Universidad Nacional
Autónoma de México IN101607.}

\maketitle

%%%%%%%%%%%%%%%%%%%%%%%%%%%%%%%%%%%%%%%%%%%%%%%%%%%%%%%%%%%%%%%%
%%%%%%%%%%%%%%%%%%%%%%%%%%%%%%%%%%%%%%%%%%%%%%%%%%%%%%%%%%%%%%%%
\section{Introduction.}

Let $\Lambda$ be an Artin algebra, and consider $\mad\Lambda$ the
class of finitely generated left $\Lambda$-modules. The finitistic
dimension of $\Lambda$ is then defined to be

$$\fd \Lambda = \sup \{\pd M \; : \; M \in\mad\Lambda \hbox{ and } \pd M <
\infty\},$$

\smallskip
\noindent where $\pd M$ denotes the projective dimension of $M$. It
was conjectured by Bass in the 60's that $\fd \Lambda$ is always
finite. Since then, this conjecture was shown to hold for many
classes of algebras \cite{AR,C,Co,IT,EHIS,Xi,Xi2, ZH0,ZH}. In particular the conjecture holds
for Artin algebras with vanishing radical cube \cite{C,GZ,ZH0}.   In this
paper, we generalize this result to Artin algebras having at most
three radical layers of infinite projective dimension.

\bigskip
\noindent{\bf Theorem.} If $_\Lambda\Lambda$ has at most three radical layers of infinite projective dimension,
then $\fd \Lambda$ is finite.

\bigskip

We also provide a bound on the finitistic dimension of $\Lambda$ under the preceding hypothesis.
In order to
achieve our goal, we use the $\Psi$ function of Igusa and Todorov \cite{IT} and introduce the notion of infinite-layer length
which counts in an efficient manner the number of (not
necessarily radical) layers of infinite projective dimension of a
module.

\section{Notations and definitions.}

For an Artin algebra $\Lambda$ and a $\Lambda$-module $M$, we denote
by $\tap M$ and $\soc M$ the top and socle of $M$ respectively. Given a
class $\mathscr C$ of objects in $\mad\Lambda$, the projective dimension of $\mathscr C$ is  $\pd \mathscr C:=\sup\{\pd M \; : \; M\in\mathscr C\}$ if the class $\mathscr C$ is not empty,
otherwise it is zero. Moreover, the {\bf{finitistic dimension}} of the class
$\mathscr C$ is $\mathrm{fin.dim.}\,\mathscr C:=\pd \{M\in\mathscr C\; : \; \pd M < \infty\}$.  We
let $\mathcal{S}^\infty$ be the finite set consisting of all isomorphism
classes of simple $\Lambda$-modules of infinite projective
dimension. Similarly, we denote by $\mathcal S^{<\infty}$ the finite set
consisting of all isomorphism classes of simple $\Lambda$-modules of
finite projective dimension. Then $\alpha:= \pd \mathcal S^{<\infty}$ is finite.
Further, we denote by $[M:\mathcal S^\infty]$ the
number of (not necessarily distinct) composition factors of $M$ of
infinite projective dimension.  Note that if all the composition
factors of $M$ belong to $\mathcal S^{<\infty}$ or $M=0$, then we
have  $[M:\mathcal S^\infty]=0$ and $\pd M \leq \alpha$.
\

 We now recall the definition and main properties of the function $\Psi$
 of Igusa and Todorov \cite{IT}.  Let $K$ denote the quotient of the free abelian group
 generated by
all the symbols $[M],$ where $M\in\mad\Lambda$, by the subgroup
generated by symbols of the form: (a) $[C]-[A]-[B]$ if $C\simeq
A\oplus B,$ and (b) $[P]$ if $P$ is projective. Then $K$ is the free
$\Bbb{Z}$-module generated by the iso-classes of indecomposable
finitely generated non-projective $\Lambda$-modules. In \cite{IT},
K. Igusa and G. Todorov define the function $\Psi:\mad\Lambda\to
\Bbb{N}$ as follows.
\

The syzygy induces a $\Bbb{Z}$-endomorphism on $K$ that will also be
denoted by $\Omega.$ That is, we have a $\Bbb{Z}$-homomorphism
$\Omega :K\rightarrow K$ where $\Omega [M]:=[\Omega\,M].$ For a
given $\Lambda$-module $M,$ we denote by $<M>$ the
$\Bbb{Z}$-submodule of $K$ generated by the indecomposable direct
summands of $M.$ Since $\Bbb{Z}$ is a Noetherian ring, Fitting's
lemma implies that there is an integer $n$ such that
$\Omega:\Omega^m<M>\;\rightarrow\Omega^{m+1}<M>$ is an isomorphism
for all $m\geq n;$ hence there exists a smallest non-negative
integer $\Phi\,(M)$ such that $\Omega:\Omega^m<M>\;\rightarrow
\Omega^{m+1}<M>$ is an isomorphism for all $m\geq \Phi\,(M).$ Let
$\mathscr{C}_M$ be the set whose elements are the direct summands of
$\Omega^{\Phi(M)}M.$ Then we set:
$$\Psi\,(M):=\Phi\,(M)+\fd\mathscr{C}_M.$$ The following result is due to K. Igusa and
G. Todorov.

\begin{pro}\cite{IT}\label{ProIT} The function $\Psi:\mad\Lambda\to\Bbb{N}$ satisfies the following properties.
\begin{itemize}
\item[(a)] If $\pd\, M$ is finite then $\Psi\,(M)=\pd\,M.$ On the other hand, $\Psi\,(M)=0$ if $M$ is indecomposable and $\pd\,M=\infty,$
\item[(b)] $\Psi\,(M)=\Psi\,(N)$ if $\add\,M=\add\,N,$
\item[(c)] $\Psi(M)\leq\Psi(M\oplus N),$
\item[(d)] $\Psi(M\oplus P)=\Psi(M)$ for any projective $\Lambda$-module $P,$
\item[(e)] If $0\to A\to B\to C\to 0$ is an exact sequence in $\mad\Lambda$ and $\pd\,C$ is finite then $\pd\,C\leq\Psi\,(A\oplus B)+1.$
\end{itemize}
\end{pro}
\

\noindent Y. Wang proved the following useful inequality, which is a
direct consequence of \ref{ProIT} (e).

\begin{lem}\label{W} \cite{W} If $0\to A\to B\to C\to 0$ is an exact sequence in $\mad\Lambda$ and $\pd\,B$ is finite then
$\pd\,B\leq 2+\Psi\,(\Omega\,A\oplus \Omega^2\,C).$
\end{lem}

We will also need the following result which appears in our previous
paper \cite{HLM}.

\begin{pro}\cite{HLM}\label{LemaH} For all $M$ in $\mad\Lambda$, $\Psi\,(M)\leq 1+ \Psi\,(\Omega\,M).$
\end{pro}

\section{The functors Q and S.}

In this section we introduce endofunctors $Q$ and $S$ on $\mad
\Lambda$ that associate to a $\Lambda$-module $M$ a quotient $Q(M)$
and a submodule $S(M)$ of $M$ with the properties that the socle of
$Q(M)$ and the top of $S(M)$
 both lie in $\add \mathcal S^\infty$ whenever $[M:\mathcal S^\infty]\neq 0$. These functors
 will be used in the definition of the infinite-layer length of a module. Throughout this section,
 we let $\alpha=\pd \mathcal S^{<\infty}$.

Given any class $\mathscr{C}$ of
$\Lambda$-modules, we will consider the full subcategory $\F(\mathscr{C})$ of  $\mathrm{mod}\,\Lambda$
whose objects are the $\mathscr{C}$-filtered
$\Lambda$-modules. That is, $M \in\F(\mathscr{C})$ if there is a finite
chain $0=M_0\subseteq M_1\subseteq \cdots\subseteq M_m=M$ of
submodules of $M$ with $m\geq 0$ and such that each quotient $M_i/M_{i-1}$ is
isomorphic to some object in $\mathscr{C}.$ For example, an object $M\in\mad\Lambda$ is filtered by $\mathcal S^{<\infty}$ if and only if $[M:\mathcal S^\infty]=0.$

\begin{lem}\label{SMAX}Let $M$ be a $\Lambda$-module.  Then the set of
all submodules of $M$ filtered by $\mathcal S^{<\infty}$ admits a
unique maximal element that will be denoted by $K(M)$.
\end{lem}
\begin{dem} Note first that this class is not empty since
$[0:\mathcal S^\infty]=0$.  Now let $K_1$ and $K_2$ be two submodules of $M$
filtered by $\mathcal S^{<\infty}$ that are maximal with this property. Since
$K_1+K_2$ is a quotient of $K_1\oplus K_2\in\mathcal{F}\,(\mathcal{S}^{<\infty})$, we get that $K_1+K_2$ is also filtered by
$\mathcal S^{<\infty}$.  Hence $K_1=K_1+K_2=K_2$.
\end{dem}

Let $Q(M)$ be the quotient $M/K(M)$, that is we have the following
exact sequence: $0\rightarrow K(M) \stackrel{i_M}{\rightarrow} M \rightarrow Q(M)
\rightarrow 0$. Given a morphism $f:M\to N$ of $\Lambda$-modules, we
have the following diagram:

$$\xymatrix{
0\ar[r] & K(M)\ar[r]^{i_M}\ar@{-->}[d]^{K(f)} & M\ar[r]\ar[d]^{f} &
Q(M)\ar[r]\ar@{-->}[d]^{Q(f)} & 0\\
0\ar[r]& K(N)\ar[r]^{i_N} & N\ar[r] & Q(N)\ar[r] & 0      }$$

\noindent Now $fi_M(K(M))$ is a quotient of $K(M)\in\mathcal F(\mathcal S^{<\infty})$ therefore it is a submodule of $N$ filtered by $\mathcal S^{<\infty}$. It then follows from the maximality of $K(N)$
that $fi_M$ factors uniquely through $i_N$, giving the map $K(f)$.  Passing to the cokernels,
we obtain a map $Q(f): Q(M) \to Q(N)$ that makes the above diagram
commute. It is then straightforward to verify that $K$ and $Q,$ as defined
above, are indeed additive functors whose main properties are listed below.
Note that $\pd K(M) \leq \alpha$.

\begin{pro}\label{QM} The functors $K,\,Q:\mad \Lambda\rightarrow \mad \Lambda$
defined above have the following properties.
\begin{itemize}
\item[(a)]$\QM=0 $ if and only if $M\in\mathcal{F}\,(\mathcal{S}^{<\infty})$.
\vspace{.2cm}

\item[(b)] If \,$\soc M \in \add \mathcal S^\infty$ then $Q(M)=M$,
\vspace{.2cm}
\item[(c)]$\pd M< \infty $ if and only if $ \pd \QM < \infty$,
\vspace{.2cm}
\item[(d)]If $[M:\mathcal S^\infty]\neq 0$ then $\soc \QM \in \add\mathcal S^\infty$,
\vspace{.2cm}
\item[(e)]$\Omega^{\alpha+1}M \oplus P \cong\Omega^{\alpha+1}\QM\oplus P'$ for some projective $\Lambda$-modules $P$ and $P'$,
\vspace{.2cm}
\item[(f)]$\pd M \leq\max\{\pd \QM, \alpha\}$,
\vspace{.2cm}
\item[(g)]If $f:M\to N$ is a monomorphism (epimorphism), then $Q(f):Q(M)\to Q(N)$ is a monomorphism (epimorphism),
\vspace{.2cm}
\item[(h)] $Q^2=Q,$ $K^2=K$ and $K\,Q=0=Q\,K.$

\end{itemize}

\end{pro}
\begin{dem} Statements (c), (e) and (f) are easily verified using the exact sequence
$0\rightarrow K(M) \rightarrow M \rightarrow \QM \rightarrow 0$ and
the fact that $\pd K(M) \leq \pd \mathcal S^{<\infty}=\alpha<\infty$.
\

(a) By definition, we have the equivalences  $0=\QM\iff K(M)=M \iff [M:\mathcal S^\infty]=0.$
\

(b) If $K(M)\neq 0$, then $0\neq \soc K(M)\subseteq \soc M\in \add S^\infty$, a
contradiction since $K(M)\in \mathcal F(\mathcal S^{<\infty})$. Thus $K(M)=0$ and $Q(M)=M$.
\

(d) Since $[M:\mathcal S^\infty]\neq 0$, we have from (a) that $\QM\neq 0$.
Assume that $\soc \QM$ admits a simple summand $S$ of finite
projective dimension, and consider the following commutative diagram

$$\xymatrix{
0\ar[r] & K(M)\ar[r]^i\ar@{=}[d] & E\ar[r]\ar[d]^{f} &
S\ar[r]\ar[d] & 0\\
0\ar[r]& K(M)\ar[r] & M\ar[r] & \QM\ar[r] & 0\\ }$$

\noindent where the upper exact sequence is obtained by pullback.
Applying the snake lemma, we infer that $f$ is a monomorphism.
Moreover, $E$ is filtered by $\mathcal S^{<\infty}$ since both $K(M)$
and $S$ are so. The maximality of
$K(M)$ then implies that $i$ is an isomorphism, thus $S=0$, a
contradiction.
\

(g) Let $f:M\to N$ be a morphism of $\Lambda$-modules. Consider the following
exact and commutative diagram
$$\xymatrix{
0\ar[r] & K(M)\ar[r]\ar[d]^{K(f)} & M\ar[r]\ar[d]^{f} &
Q(M)\ar[r]\ar[d]^{Q(f)} & 0\\
0\ar[r]& K(N)\ar[r] & N\ar[r] & Q(N)\ar[r] & 0\\ }$$ If $f:M\to N$ is an epimorphism, then from the
diagram above, we see that $Q(f)$ is also an epimorphism.\\
Suppose that $f:M\to N$ is a monomorphism.  If $M\in\mathcal{F}(\mathcal S^{<\infty})$,
then by (a) $Q(M)=0;$ and hence $Q(f)$ is a
monomorphism. Assume now that $[M:\mathcal S^\infty]\neq 0$. It then follows from (d) that $\soc Q(M) \in
\add \mathcal S^\infty$. Applying the snake lemma to the diagram above, we get a
monomorphism from $X:=\Ker Q(f)$ to $Y:=\Coker K(f)$. If $X\neq 0$, then
$0\neq \soc X \subseteq  \soc Q(M) \in \add \mathcal S^\infty$.
Also, $\soc X \subseteq \soc Y$, a contradiction since $Y\in\mathcal{F}(\mathcal S^{<\infty})$. So $X=0$ and $Q(f)$ is a monomorphism.
\

(h) By taking the pull-back to the canonical morphisms $M\to Q(M)\leftarrow K\,Q(M),$ we get the following exact and commutative diagram

$$\xymatrix{
&  & 0 \ar[d] & 0 \ar[d] \\
0\ar[r] & K(M)\ar[r]^j \ar@{=}[d]  & E \ar[r]\ar[d] & K\,Q(M)\ar[r]\ar[d] & 0\\
0\ar[r]& K(M)\ar[r]  & M \ar[r]\ar[d]  & Q(M)\ar[r] \ar[d] & 0 \\
&  & Q^2(M)\ar@{=}[r] \ar[d] & Q^2(M) \ar[d]\\
&  & 0 & 0     }$$

Hence, by \ref{SMAX} and the fact that $E\in\mathcal{F}(\mathcal S^{<\infty}),$ we conclude that $j$ is an isomorphism so that $K\,Q(M)=0$ and $Q(M)=Q^2(M).$
On the other hand, since $K(M)\in\mathcal{F}(\mathcal S^{<\infty}),$ we get from \ref{QM}(a) that $Q\,K(M)=0$ and $K^2(M)=K(M)$.
\end{dem}

We now proceed to give the construction of the functor $S$ which is
dual to $Q$.  We omit the proofs since they are essentially the same
as those we previously did. In what follows, a quotient of $M$ is an epimorphism $g:M\to C;$ moreover, if $g_1:M\to C_1$ and $g_2:M\to C_2$ are two quotients of $M$, we say that $g_1$ is {\em greater than or equal to} $g_2$ if there is an epimorphism $h:C_1\to C_2$ such that $hg_1=g_2.$ Finally, we say that a quotient $g:M\to C$ is filtered by $\mathcal S^{<\infty}$ in case $C$ is so.

\begin{lem}\label{QMAX}Let $M$ be a $\Lambda$-module.  Then the set of
all quotients of $M$ filtered by $\mathcal S^{<\infty}$ admits a
unique (up to isomorphism) maximal element that will be denoted by $p_M:M\to C(M)$.
\end{lem}

Consider the exact sequence $0\rightarrow \SM\rightarrow M \stackrel{p_M}{\rightarrow}
C(M) \rightarrow 0$. Given a morphism of
$\Lambda$-modules $f:M\to N$, we have the following diagram:

$$\xymatrix{
0\ar[r] & S(M)\ar[r]\ar@{-->}[d]^{S(f)} & M\ar[r]^{p_M}\ar[d]^{f} &
C(M)\ar[r]\ar@{-->}[d]^{C(f)} & 0\\
0\ar[r]& S(N)\ar[r] & N\ar[r]^{p_N} & C(N)\ar[r] & 0\\ }$$

\noindent Now $p_Nf(M)$ is a submodule of
$C(N)\in\mathcal{F}(\mathcal S^{<\infty})$, therefore $M\to p_Nf(M)$
is a quotient of $M$ filtered by $\mathcal S^{<\infty}$. It then
follows from the maximality of $p_M:M\to C_M$ that $p_Nf$ factors
uniquely through $p_M$ giving us the map $C(f)$. By passing to
kernels, we obtain the map $S(f)$ that makes the above diagram
commute. It is now straightforward to verify that $S$ and $C$ are
additive functors.  We then have the following properties of such functors,
the proof of which is dual to \ref{QM}.

\begin{pro}\label{SM}The functors $S,$ $C:\mad \Lambda\to \mad\Lambda$,
defined above, have the following properties.
\begin{itemize}
\item[(a)]$\SM=0 $ if and only if $M\in\mathcal{F}(\mathcal S^{<\infty})$,
\vspace{.2cm}
\item[(b)] If $\tap M\in\add \mathcal S^\infty$, then $S(M)=M$,
\vspace{.2cm}
\item[(c)]$\pd M< \infty $ if and only if $ \pd \SM < \infty$,
\vspace{.2cm}
\item[(d)]If $[M:\mathcal S^\infty]\neq 0$ then $\tap\SM \in \add\mathcal S^\infty$,
\vspace{.2cm}
\item[(e)]$\Omega^{\alpha}M \oplus P \cong\Omega^{\alpha}\SM\oplus P'$ for some projective $\Lambda$-modules $P$ and $P'$,
\vspace{.2cm}
\item[(f)]$\pd M \leq\max\{\pd \SM, \alpha\}$,
\vspace{.2cm}
\item[(g)]If $f:M\to N$ is a monomorphism (epimorphism), then $S(f):S(M)\to S(N)$ is a monomorphism (epimorphism),
\vspace{.2cm}
\item[(h)] $C^2=C,$ $S^2=S$ and $C\,S=0=S\,C$.
\end{itemize}
\end{pro}
\

Note that it can also be shown that $QS$ and $SQ$ are naturally isomorphic functors.

\section{The infinite-layer length of a module.}

In this section, we introduce the invariant $\LL^\infty(M)$ for a
$\Lambda$-module $M$.  Our goal is to count the number of "layers"
of infinite projective dimension of $M$. A natural way to proceed
would be to consider radical layers.  Recall that a {\bf{radical layer}} of
a module $M$ is a semisimple module of the form
$\rad^iM/\rad^{i+1}M$ for some $0\leq i <\ell(M)$, where $\ell(M)$ is
the Loewy length of $M$. It then seems reasonable to define the
"infinite-layer length" $\ell^\infty(M)$ of a module $M$ as follows.

\begin{defi} For any $M\in\mathrm{mod}\,\Lambda,$ we set
$\ell^\infty(M)$ to be the number of radical layers of $M$ that have
infinite projective dimension.
\end{defi}
 As natural as it seems, this
definition has some flaws. For example, if $K$ is a submodule of $M$, we do not always
have that $\ell^\infty(K) \leq \ell^\infty(M)$.

We will now define another "infinite-layer length", which we will denote by $\LL^\infty$, that satisfies
the above stated property and is better than $\ell^\infty$ in the
sense that $\LL^\infty(M) \leq \ell^\infty(M)$. As an analogy, the
Loewey length of a module $M$ can be defined to be the smallest
nonnegative integer $i$ such that $\rad^iM=0$.  In our case, since
we are only interested in layers of infinite projective dimension,
we use the functor $S$ to "level" our module prior to taking the
radical.  This guarantees that at each step, a layer consisting of a maximal number
of simples of infinite projective dimension is removed.

Given a module $M$, we start by calculating $S(M)=M^0$. Unless
$S(M)=0$, the top of $M^0$ lies in $\add \mathcal S^\infty$. We
then take the radical of $M^0$, and let $M^1:= S(\rad(M^0))$.
Iterating the process, we let $M^{i+1}:=S(\rad(M^i))$ for all $i\geq
0$. This procedure leads us to consider the following additive
functor
$$F:=S\circ \rad:\mathrm{mod}\,\Lambda\to\mathrm{mod}\,\Lambda,$$ where $S$ is
the functor defined in the previous section and $\rad$ is the
radical functor. Note that $F^0$ is the identity functor; and so, by
using the functor $F,$ it is now easy to see that $M^i=F^i\,(S(M))$ for all $i\geq 0$.

\begin{defi}\label{ll} Let $\Lambda$ be an Artin algebra.  For a finitely
generated $\Lambda$-module $M,$ we define the {\bf\em infinite-layer length} of $M$ to be
$$\LL^\infty(M):=\mathrm{min}\,\{i\geq 0\; : \; F^i (S(M))=0\}.$$
\end{defi}

Note that it follows directly from the definition that for all $M\in\mad \Lambda$, $\LL^\infty(M) \leq \ell^\infty(M)$.

\begin{ex}\label{ex}
Consider the bound quiver algebra $\Lambda=kQ/I$ where $Q$ is given by

$$
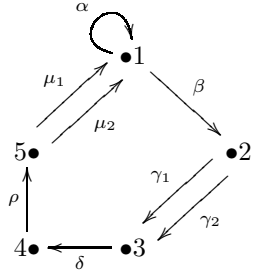
\xymatrix{ & {\bullet}1\ar@(l,u)[]^{\alpha} \ar[dr]^{\beta} & \\
{5}{\bullet}\ar@<1ex>[ru]^{\mu_1} \ar@<-1ex>[ru]_{\mu_2}& &{\bullet}2 \ar@<1ex>[dl]^{\gamma_2} \ar@<-1ex>[dl]_{\gamma_1}\\
{4}{\bullet}\ar[u]^{\rho} & {\bullet}{3} \ar[l]^{\delta} & }$$

\noindent and $I$ is generated by the set of paths $\{\alpha^3 , \alpha\beta, \rho\mu_i\alpha,\mu_i\beta ,\gamma_1\delta-\gamma_2\delta\}$ where $i=1,2$. Then $\mathcal S^{\infty}=\{S(1), S(4)\}$ and the indecomposable projective $\Lambda$-modules are

$$ 
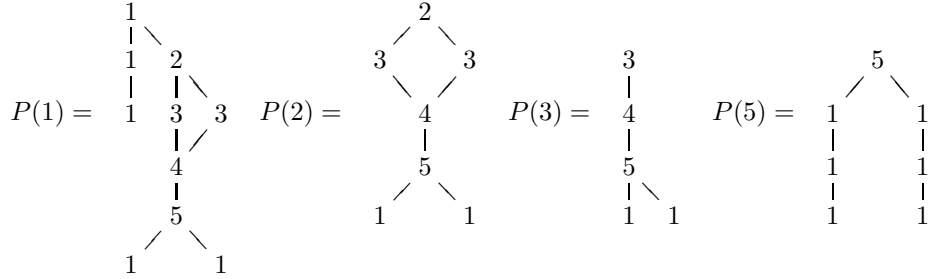
\xymatrix@-1.5pc {&1\edge[d]\edge[dr]&& &&& 2\edge[dl]\edge[dr]&&&&&&&&\\
                     &1\edge[d]&2\edge[d]\edge[dr]&&&3\edge[dr]&&3\edge[dl]&      & 3\edge[d]&&&&5\edge[dl]\edge[dr]\\ P(1)=&1&3\edge[d]&3\edge[dl]&P(2)=&&4\edge[d]&&P(3)=&4\edge[d]&&P(5)=&1\edge[d]&& 1\edge[d]\\
                     &&4\edge[d] &&&& 5\edge[dl]\edge[dr] &&&5\edge[d]\edge[dr]&&&1\edge[d]&&1\edge[d]\\
                     && 5\edge[dl] \edge[dr] &           && 1 &  &1&         &1&1&&1&&1\\
                     &1 & & 1}$$

\noindent and $P(4)=\rad P(3)$.  In order to compute $\LL^\infty(P(1))$, we need to find the smallest nonnegative integer $i$ such that $F^i(S(P(1)))=0$. Note that $S(P(1))=P(1)$ since $\tap P(1)=S(1)\in \add \mathcal S^\infty$. Thus $FS(P(1))=F(P(1))=S \rad (P(1))$ is the following module:
$$\xymatrix@-1.5pc{&& && 4\edge[d]\\
                   FS(P(1))=&&1\edge[d]&\oplus &5\edge[dr]\edge[dl]&\\
                   &&1&1&&1}$$

\noindent Continuing this process, we obtain $F^2S(P(1))= (S(1))^3$ and $F^3S(P(1))=0$, thus $\LL^\infty(P(1))=3$.  It is easy to see that among the six radical layers of $P(1)$ only the fifth has finite projective dimension,
therefore $\ell^\infty(P(1))=5$.  Further computations show that $\LL^\infty(P(2))=\LL^\infty(P(3))=\LL^\infty(P(4))=2$ and that $\LL^\infty(P(5))=3$.
\end{ex}

\

\begin{lem}\label{LemaF} Consider the functor $F^i:\mad\Lambda\to\mad\Lambda$ for $i\geq 0.$ Then
 $F^i$ and $F^i\circ S$ preserve monomorphisms and epimorphisms.
\end{lem}
\begin{dem} From \ref{SM} (g), we know that $S$ preserves both
monomorphisms and epimorphisms.  This is also true of the functor
$\rad$.  Consequently, for all non negative integers $i$, we have
that $F^i$ and $F^i\circ S$ also preserves monomorphisms and
epimorphisms.
\end{dem}

The next result shows that $\LL^\infty$ behaves naturally with monomorphisms, epimorphisms and direct sums.

\begin{pro}\label{ll1}  Let $\Lambda$ be an Artin algebra and $L,\, M$ and $N$ be
$\Lambda$-modules.
\begin{itemize}
\item[(a)] If $f:L\to M$ is a monomorphism, then $\LL^\infty(L)\leq
\LL^\infty(M)$,
\vspace{.2cm}
\item[(b)] If $g:M\to N$ is an epimorphism, then $\LL^\infty(N)\leq
\LL^\infty(M)$,
\vspace{.2cm}
\item[(c)] $\LL^\infty(M\oplus N)=\mathrm{max}\,(\LL^\infty(M),\,\LL^\infty(N)),$
\vspace{.2cm}
\item[(d)] $\LL^\infty(M)\leq\LL^\infty({}_\Lambda\Lambda).$

\end{itemize}
\end{pro}
\begin{dem} (a) Let $f:L\to M$ be a monomorphism and assume that $\LL^\infty(M)=m$. Then
$F^m(S(M))=0;$ and so by \ref{LemaF} we get
$F^m(S(L))=0.$ Thus $\LL^\infty(L)\leq m$, proving (a). A similar
argument holds for (b).
\

(c) Since the functor $F^i\,S$ is aditive, we have that $$F^i\,S(M\oplus N)\simeq F^i\,S(M)\oplus F^i\,S(N),$$ therefore $\mathrm{max}\,(\LL^\infty(M),\,\LL^\infty(N))\leq\LL^\infty(M\oplus N).$ The reverse inequality follows
from the fact  that $F^{i_0}S(X)=0$ implies $F^iS(X)=0$ for all $i\geq i_0.$
\

(d) Using that $M$ is a quotient of ${}{}_\Lambda\Lambda^n$ for some natural number $n,$ we obtain from (b) that $\LL^\infty(M)\leq\LL^\infty({}_\Lambda\Lambda^n);$ thus (d) follows from (c).
\end{dem}

\begin{pro}\label{ll2}
Let $\Lambda$ be an Artin algebra and $M$ a $\Lambda$-module.
\begin{itemize}
\item[(a)] $\LL^\infty(M)=0$ if and only if $M\in\mathcal{F}(\mathcal S^{<\infty})$,
\vspace{.2cm}
\item[(b)] $\LL^\infty(S(M))=\LL^\infty(M)$,
\vspace{.2cm}
\item[(c)]If $M\neq 0$ and $\tap M$ lies in $\add \mathcal S^\infty$ then $\LL^\infty
(\rad M) =\LL^\infty(M)-1$.
\end{itemize}
\end{pro}
\begin{dem} (a) Let $M\in\mathcal{F}(\mathcal S^{<\infty}).$ Then by \ref{SM}(a) we have $S(M)=0$, thus $F^0(S(M))=0$ and $\LL^\infty(M)=0$.  Conversely, if $\LL^\infty(M)=0$ then $S(M)=0$ and hence $M\in\mathcal{F}(\mathcal S^{<\infty}).$
\

(b) This follows from the definition and the fact that $S^2(M)=S(M)$.
\

(c) It follows from (a) that $\LL^\infty(M)>0$ and from \ref{SM}(b)
that $S(M)=M$. For all $i\geq 1$, we have
$F^i(S(M))=F^i(M)=(S\circ\rad)^{i-1}S\,(\rad M)=F^{i-1}(S(\rad
M)).$ Therefore $\LL^\infty(\rad M)=\LL^\infty(M)-1$.
\end{dem}

%Redundant, this is just a reformulation of ll3 and ll4
%\begin{cor}\label{coroll4} Let $\Lambda$ be an Artin algebra and $M$ be a finitely generated $\Lambda$-module. If %$\LL^\infty(M)>0,$ then
%\begin{itemize}
%\item[(a)] $\LL^\infty (G(M))=\LL^\infty(Q(M))-1,$
%\vspace{.2cm}
%\item[(b)] $\LL_\infty (G(M))=\LL_\infty(M)-1.$
%\end{itemize}
%\end{cor}
%\begin{dem} Since $\LL^\infty(M)>0,$ we get, from \ref{ll2} (a), that $[M:\mathcal S^\infty]\neq 0;$ and so, by %\ref{QM} (d), we conclude that $\soc Q(M)\in \add \mathcal S^\infty]$ Hence, (a) follows from \ref{ll3}; and moreover, %(b) is a consequence of \ref{ll4} (a) and \ref{ll4} (b).
%\end{dem}
%\

\section{Main results.}

We now proceed to show that an Artin algebra $\Lambda$ with infinite-layer length at most three has finite
finitistic dimension. Throughout this section, we will use the notation $\alpha=\pd \mathcal S^{<\infty}$ and $\Sigma=\oplus_{S\in\mathcal S^\infty}S.$

\begin{defi} Let $\Lambda$ be an Artin algebra and $M$ a finitely generated $\Lambda$-module such that $[M:\mathcal S^\infty]\neq 0$. Then, we set
$$S^\infty_M:=\bigoplus_{S\in\mathcal S^\infty}\,S^{[M:S]}\in\add \mathcal S^\infty.$$
\end{defi}

\begin{lem}\label{ONE}  Let $\Lambda$ be an Artin algebra and $M$ a finitely generated $\Lambda$-module. If $\LL^\infty(M)=1,$ then
\begin{itemize}
\item[(a)] $\pd M=\infty,$ and $S^\infty_M={S}_{S(M)}^{\infty}=\tap S(M)$,
\vspace{.2cm}
\item[(b)] $\Omega^{\alpha+1}(M)\oplus P\simeq \Omega^{\alpha+1}( S^\infty_M)\oplus P'$ for some projective $\Lambda$-modules $P$ and $P'.$
\end{itemize}
\end{lem}
\begin{dem} Suppose that $\LL^\infty(M)=1.$ Then \ref{ll2} (b) and (c) give $\LL^\infty (\rad S(M))=0;$ and so, by \ref{ll2}(a), we conclude that $\rad S(M) \in\F(\mathcal S^{<\infty}).$ Therefore, by considering the following exact sequence
$$0\rightarrow \rad S(M) \rightarrow S(M) \rightarrow \tap S(M) \rightarrow 0,$$ we obtain that $\tap S(M)$ contains all composition
 factors of $S(M)$ (and hence of $M$) of infinite projective dimension. This proves  the second part of (a) since $\tap S(M)\in\add \mathcal S^{\infty}$ (see \ref{SM}(d)).
\

On the other hand, since $\pd \rad S(M) \leq \alpha$ and $\pd \tap S(M)=\infty$, we have $\pd S(M)=\infty$ and hence from \ref{SM}(c), $\pd M= \infty$.  Also, the above exact sequence yields $\Omega^{\alpha+1}
(S(M))\oplus P\simeq \Omega^{\alpha+1} (\tap S(M)) \oplus P'$ for some projective $\Lambda$-modules $P$ and $P'.$ Thus (b) follow from \ref{SM}(e).
\end{dem}

We now consider the case when $\LL^\infty(M)=2$.

\begin{pro}\label{TeoOne} Let $\Lambda$ be an Artin algebra. If $M$ is a module of finite projective dimension and $\LL^\infty(M)=2$, then $\pd M \leq \alpha +2 +\Psi(\Omega^{\alpha+1}(\Sigma)\oplus
\Omega^{\alpha+2}(\Sigma)).$
\end{pro}

\begin{dem} Let $M$ be of finite projective dimension and of infinite-layer length $2$.  Then by \ref{SM}(c), $\pd S(M)$ is finite. Applying \ref{ll2}(b) and (c), we obtain $\LL^\infty(\rad S(M))=1.$ Hence, from \ref{ONE}(b), we get
$$\Omega^{\alpha+1}(\rad S(M))\oplus P\simeq \Omega^{\alpha+1}( S^\infty_{\rad S(M)})\oplus P'$$ for some projective $\Lambda$-modules $P$ and $P'.$ Applying \ref{W} and \ref{LemaH} to the exact sequence
$0\rightarrow \rad S(M) \rightarrow S(M) \rightarrow \tap S(M) \rightarrow 0,$ we get
\begin{eqnarray*}
\pd S(M) &\leq& 2 + \Psi[\Omega (\rad S(M)) \oplus \Omega^2 (\tap S(M))]\\
&\leq & 2+\alpha +
\Psi[\Omega^{\alpha+1} (\rad S(M)) \oplus \Omega^{\alpha+2} (\tap S(M))]\\
&=&2+\alpha + \Psi[\Omega^{\alpha+1} (S^\infty_{\rad S(M)}) \oplus
\Omega^{\alpha+2} (\tap S(M))]\\
&\leq& 2+\alpha + \Psi[\Omega^{\alpha+1}(\Sigma) \oplus
\Omega^{\alpha+2}(\Sigma)],
\end{eqnarray*}
\noindent where the last inequality follows from \ref{ProIT} (c)
and (b).  It then follows from \ref{SM} (f) that $$\pd M
\leq \max\{\pd S(M), \alpha\} \leq 2+\alpha + \Psi(\Omega^{\alpha+1}
(\Sigma) \oplus \Omega^{\alpha+2} (\Sigma)),$$ proving the result.
\end{dem}
\

In order to treat the case when $\LL^\infty(M)=3$, we will need to following lemma.

\begin{lem}\label{tres} Let $\Lambda$ be an Artin algebra such that $\mathcal S^\infty\neq \emptyset$, and $M\in\mad\Lambda$,  then
$$\LL^\infty(\Omega S(M))\leq \LL^\infty(_\Lambda\Lambda)-1.$$
\end{lem}
\begin{dem} Let $M\in\mad\Lambda.$ If $[M:\mathcal S^\infty]=0,$ then we have by \ref{SM}(a) that $\Omega S(M)=0$, and $\LL^\infty (0)=0 \leq \LL^\infty(_\Lambda\Lambda)-1$ since $\mathcal S^\infty\neq\emptyset$.

Suppose that $[M:\mathcal S^\infty]\neq 0$ and let $P$ be the projective cover of $S(M).$ On one hand, we have from \ref{ll1}(d) that $\LL^\infty(P)-1\leq \LL^{\infty} (_\Lambda\Lambda)-1$. On the other hand, using \ref{SM}(d), we get that $\tap P\simeq\tap\SM\in\add \mathcal S^\infty;$ thus, by \ref{ll2}(c), we have $\LL^\infty(\rad P)=\LL^\infty(P)-1$. Finally, since $\Omega S(M)\subseteq \rad P,$ we obtain using \ref{ll1}(a) that $\LL^\infty(\Omega S(M))\leq\LL^\infty(\rad P)$ proving the result.
\end{dem}
\

We are now in position to prove our main result.

\begin{teo}\label{main} Let $\Lambda$ be an Artin algebra such that $\LL^\infty(_\Lambda \Lambda)\leq 3$, then
$$\fd\Lambda\leq 3+\alpha + \Psi(\Omega^{\alpha+1}
(\Sigma) \oplus \Omega^{\alpha+2} (\Sigma)),$$ where $\alpha=\pd \mathcal S^{<\infty}$  and $\Sigma=\oplus_{S\in\mathcal S^\infty}S.$
\end{teo}
\begin{dem} Let $M\in\mad\Lambda$ be of finite projective dimension. It follows from \ref{ll1}(d) that $\LL^\infty(M)\leq 3$ and from \ref{ONE} that $\LL^\infty(M)\neq 1$.  If $\LL^\infty(M)=0$, then $\pd M\leq \alpha$.
If $\LL^\infty(M)=2$, then it follows from \ref{TeoOne} that $\pd M\leq \alpha +2 +\Psi(\Omega^{\alpha+1}(\Sigma)\oplus
\Omega^{\alpha+2}(\Sigma)).$  If $\LL^\infty(M)=3$, then $\pd M\leq \max\{\alpha, \pd S(M)\}\leq \max\{\alpha, 1+ \pd \Omega(S(M))\}\leq 3+\alpha + \Psi(\Omega^{\alpha+1}
(\Sigma) \oplus \Omega^{\alpha+2} (\Sigma))$ since by \ref{tres}, $\LL^\infty(\Omega(S(M))\leq 2$.
\end{dem}

\begin{ex} Consider the algebra $\Lambda$ of example \ref{ex}. We showed in this example that $$\max\{\LL^\infty(P(i)) \, :\, 1\leq i\leq 5\}=3,$$ therefore we infer from \ref{ll1}(c) that $\LL^\infty(_\Lambda\Lambda)= 3$.  Also, $\mathcal S^\infty=\{S(1),S(4)\}$ and $\alpha= \pd\mathcal S^{<\infty}=2$. Using \ref{ProIT} we get $\Psi[\Omega^{3}(S(1)\oplus S(4))\oplus \Omega^{4}(S(1)\oplus S(4))]=\Psi(S(1)\oplus T)=0$, where $T$ is the two-dimensional indecomposable module whose top and socle are both isomorphic to $S(1)$.  It then follows from \ref{main} that $\fd\Lambda \leq 5$.  Note that this algebra has 5 radical layers of infinite projective dimension and so $\ell^\infty (\Lambda)=5$.
\end{ex}

\
Since $\LL^\infty(M)\leq \ell^\infty(M)$ for all $M$ we immediately get the corollary.

\begin{cor} If $\Lambda$ is an Artin algebra such that
${}_\Lambda\Lambda$ has at most three radical layers of infinite projective dimension, then $\fd \Lambda$ is finite.
\end{cor}

\end{document}